\documentclass[12 pt]{article}
\usepackage{amssymb}
\usepackage[centerta gs]{amsmath}
\usepackage{amscd, amsthm}

\newtheorem{fact}{fact}[section]

\newtheorem{thm}[fact]{Theorem}
\newtheorem{lemma}[fact]{Lemma}

\newtheorem{corollario}[fact]{Corollary}
\newtheorem{defini}[fact]{Definition}
\newtheorem{osserva}[fact]{Remark}

\begin{document}
\title{Comparing $\mathbb C$ and Zilber's  Exponential Fields:  Zero Sets of Exponential Polynomials}
\author{P. D'Aquino\footnote{Department of Mathematics, Seconda Universit\`{a} di Napoli, Via Vivaldi
43, 81100 Caserta, paola.daquino@unina2.it}, A. Macintyre
\footnote{School of Mathematical Sciences, Queen Mary University
of London, Mile End Road, London E1 4NS, A.Macintyre@qmul.ac.uk} and
G. Terzo\footnote{Department of Mathematics, Seconda
Universit\`{a} di Napoli, Via Vivaldi 43, 81100 Caserta,
 giuseppina.terzo@unina2.it}}

\maketitle

\begin{abstract}
We continue the research programme of comparing the complex exponential with  Zilber's exponential.  For the latter we prove, using diophantine geometry, various properties about zero sets of exponential functions, proved for $\mathbb C$ using analytic function theory, e.g.  the Identity Theorem. 
\end{abstract}


\noindent AMS 2010  {\it Mathematics subject classification:} $03$C$60$, $11$U$09$.

\section{Introduction}

In \cite{zilber} Zilber introduced, and studied deeply,  a class of exponential fields now known as Zilber fields. There are many novelties in his analysis, including a reinterpretation of Schanuel's Conjecture in terms of Hrushovski's very general theory of predimension and strong extensions. By now there is no need to spell out yet again all his ingredients and results (see \cite{zilber}, \cite{pag}). The most dramatic aspect is that his fields satisfy Schanuel's Conjecture and  a beautiful {\it Nullstellensatz} for exponential equations. Moreover, in each uncountable cardinal there is a privileged such field, satisfying a {\it countable closure condition} and a stregthened Nullstellensatz. {\it Privileged} means that the structure in each uncountable cardinal is unique up to isomorphism. The one in cardinal $2^{\aleph_0}$ is called $\mathbb B$ by us. Zilber conjectured that  $\mathbb B \cong \mathbb C$ as exponential fields. This would, of course, imply that $\mathbb C$ satisfies Schanuel's Conjecture, and Zilber's Nullstellensatz, which seem far out of reach of current analysis of several complex variables.
\par
Zilber fields are constructed model-theoretically and have no visible topology except an obvious exponential Zariski topology. Zilber's countable closure condition is somehow an analogue of the separability of the complex field, and the fact that any finite system of exponential equations has only countably many isolated points.  Isolation in $\mathbb C$ relates, via the Jacobian criterion, to definability and dimension in $\mathbb B$, see \cite{zilber}.  
\par
We have undertaken a research programme of taking results from $\mathbb C$, proved using analysis and/or topology, and seeking exponential-algebraic proofs in $\mathbb B$. An early success was a proof in $\mathbb B$ of the Schanuel's Nullstellensatz \cite{pag}, proved in $\mathbb C$ (\cite{henrub}) using  Nevanlinna theory. More recently in \cite{DMT} we derived, in an exponential-algebraic way, Shapiro's Conjecture from Schanuel's Conjecture, thereby getting Shapiro's Conjecture in $\mathbb B$.  In the present paper we put the ideas of \cite{DMT}  to work on some problems connected to the {\bf Identity Theorem} of complex analysis \cite{alhfors}. That fundamental theorem says in particular that if the zero set of an entire function $f$ has an accumulation point then $f\equiv0$. We specialize to exponential functions, and face the obvious difficulty that ``accumulation points" has no general meaning in Zilber fields. 
\par
We prove, not only for Zilber fields, but for the much more general classes of $LEC$-fields and $LECS$-fields (see Section 2) various results proved for $\mathbb C$ using the Identity Theorem. Our replacement techniques come from diophantine geometry. Our results are not confined to ones proved in $\mathbb C$ using the Identity Theorem, see for example Theorem \ref{cyclicgroupcase}. This theorem is not true for all entire functions, see Section \ref{rationalsolutions}. An analysis of the location of the zeros   of exponential polynomials in the complex plane was obtained  by P\'olya et al. (see e.g.  \cite{Polya}) in terms of lines in the plane determined by the polynomial itself. For more recent work on this see \cite{vanderpoortentijdeman}.  Notice that this analysis makes no  sense in Zilber fields, as the very notion of a line (defined over $\mathbb R$) makes no sense.


\section{Exponential fields and exponential polynomials}
We will be working over an algebraically closed field $K$ of characteristic $0$, with a surjective exponential map whose kernel is an infinite cyclic group which will be denoted by the standard notation $2\pi i\mathbb Z$.  Here $\pi$ has a well defined meaning, see  \cite{kmo}. We will call these fields $LEC$-fields. In some cases we will assume that the field $K$ satisfies the following transcription of Schanuel's Conjecture from transcendental number  theory: 

\noindent ({\bf SC})  If $\alpha_1, \ldots , \alpha_n \in K$   then $$ tr.d._{\mathbb Q}(\alpha_1, \ldots , \alpha_n, e^{\alpha_1}, \ldots , e^{\alpha_n })\geq l.d._{\mathbb Q}(\alpha_1, \ldots , \alpha_n). $$ 
We will refer to these fields  as $LECS$-fields. The class of $LECS$-fields includes the exponential fields  introduced by Zilber in \cite{zilber}. 
\par
We will consider exponential polynomial functions over $K$ of the following form
\begin{equation}
\label{poly}
f(z) = \lambda_1e^{\mu_1 z} + \ldots + \lambda_Ne^{\mu_N z}, \mbox{ where } \lambda_i, \mu_i \in K.
\end{equation} 

The set of these polynomials form a ring $\cal{E}$ under the usual addition and multiplication. We will study subsets of the zero set of polynomials in $\mathcal E$. We will denote the zero set of $f$ by $Z(f)$.

\par
We recall some basic definitions and results for the exponential polynomials in the ring  $\cal{E}.$ The units in  $\cal{E}$ are the purely exponential terms $e^{\mu z}$ where $\mu \in K.$

\begin{defini}

An element $f$ in  $\cal{E}$ is irreducible, if there are no non-units $g$ and $h$ such that $f = gh.$
\end{defini}

\begin{defini}
Let $f =\sum_{i=1}^{N} \lambda_{i} e^{\mu_iz}$ be an exponential polynomial. The support of $f,$ denoted by $supp(f),$ is the $\Bbb Q$-space generated by $\mu_1, \ldots, \mu_N.$
\end{defini}

\begin{defini} An exponential polynomial $f (z)$ in  $\cal{E}$  is simple if $\dim(supp(f )) = 1.$ \end{defini}

It is easily seen that, up to a unit, a simple exponential polynomial is a polynomial in $e^{\mu z}$, for some $\mu \in K.$ An important example of a simple exponential polynomial is  $\sin(\pi z).$ A simple polynomial  $f$ can be factored,  up to a unit in $\mathcal E$, as $f=\prod (1-ae^{\mu z})$, where $a, \mu \in K$. If a simple polynomial $f$ has infinitely many roots then by the Pigeon-hole Principle  one factor of $f$, say $1-ae^{\mu z}$ has infinitely many zeros and these are of the form $z=(2k\pi i-\log
a)/\mu$ with $k\in \mathbb Z$,  for a fixed value of $\log a$. 
\par
We will refer to $\lambda_1e^{\mu_1z}+\lambda_2e^{\mu_2z}$ (or to the equivalent form  $1-ae^{\mu z}$)  as a simple polynomial of length $2$.  We have a complete description of 
the zero set of these polynomials.
\begin{lemma}
\label{simplepoly}
If $f(z)=\lambda_1e^{\mu_1z}+\lambda_2e^{\mu_2z}$ then 
 $Z(f)$ has dimension less or equal than $2$ over $\mathbb Q$. Moreover,  $Z(f)=A_{f} + \mathbb ZB_f$, where $A_f=\frac{\log(-\lambda_1^{-1}\lambda_2)}{\mu_1-\mu_2}$ and $B_f=\frac{2i\pi}{\mu_1-\mu_2}$. So, the zero set of $f$  is a translate of a rank  $1$ free abelian group.   
\end{lemma}

\par 
It seems that the first to consider a factorization theory for exponential polynomials over $\mathbb C$ was Ritt in \cite{ritt}.  His original idea was to reduce the factorization of an exponential polynomial to that of a classical polynomial in many variables by replacing the variables with their powers. Let $f(z) = \lambda_1e^{\mu_1 z} + \ldots + \lambda_Ne^{\mu_N z}$ be an exponential polynomial where $\lambda_i, \mu_i \in \Bbb C,$ and let $b_1, \ldots, b_D$ be a basis of the $\Bbb Z$-module spanned by  $\mu_1, \ldots ,\mu_N$. Let $Y_i = e^{b_iz},$ with $i = 1, \ldots, D.$ If each $\mu_i$ is expressed in terms of the $b_i'$s, then $f(z)$ is transformed into a classical Laurent polynomial $F(Y_1, \ldots, Y_D) \in \Bbb Q(\overline{\lambda})[Y_1, \ldots, Y_D].$ 
Clearly any factorization of $f$ produces a factorization of $F(Y_1, \ldots, Y_D).$ 

\par In general, an irreducible classical polynomial $F(Y_1, \ldots, Y_D)$ can become reducible after a substitution  of the variables by powers. 

\begin{defini} A polynomial $F(Y_1, \ldots, Y_D)$ is power irreducible  if for each $n_1, \ldots, n_D$  in $ \mathbb N - \{0\}$,  $F(Y_1^{n_1}, \ldots, Y_D^{n_D})$ is irreducible. \end{defini}

Ritt saw the importance of understanding the ways in which an irreducible classical polynomial $F(Y_1, \ldots, Y_D)$ can become reducible when the variables are replaced by their powers. His analysis gave the following result.

\begin{thm}
\label{ritt}
Let $f(z) = \lambda_1e^{\mu_1 z} + \ldots + \lambda_Ne^{\mu_N z}$  where $\lambda_i, \mu_i \in \Bbb C.$ Then $f(z)$ can be written uniquely up to order and multiplication by units as 
$$f(z) = S_1 \cdot \ldots \cdot S_k \cdot I_1 \cdot \ldots \cdot I_m,$$ where $S_j$ are simple polynomials with different supports and the $I_h$ are irreducible in $\mathcal E$.
\end{thm}

In our setting we use an analogous  result for exponential polynomials as in (\ref{poly}) over  any algebraically closed field of characteristic 0 carrying an exponentiation (see also \cite{paolagiusy}). Notice that the factorization theorem of Ritt is a result on the free $E$-ring over exponential fields, and need  not involve any analysis of the zero set of an exponential polynomial.

\section{Zero sets in $\mathbb C$}

We are interested in putting restrictions on infinite subsets of the zeros sets of  certain totally defined functions, like exponential functions.  These satisfy special properties such as Schanuel's Nullstellensatz (see \cite{pag} and \cite{henrub}). For such functions in $\mathbb C$ one may use the Identity Theorem to get information about the zero  set.  We show below that in some cases that we list, and related ones, the topology on $\mathbb C$ is not needed and the use of Identity Theorem can be replaced by diophantine geometry. Here are the examples we will consider:
\begin{enumerate}
\item 
Let $X\subseteq \mathbb Q$. There is a unique copy of the  field of rationals  in $K$ since it is a  characteristic $0$ field.  We will say that a subset $X$ of the rationals {\it accumulates}  if there exists  a Cauchy sequence of distinct elements of $X$. Note that this makes sense since the definition is given inside $\mathbb Q$, and in $\mathbb Q$ there is a perfectly good  notion of a Cauchy sequence. $\mathbb Q$  clearly accumulates in this sense, while $\mathbb Z$ and $\frac{1}{N}\mathbb Z$, where $N\in \mathbb N-\{ 0\}$   do not. Note that we are not explicitly considering the question whether a subset $X$ of $\mathbb Q$ that accumulates has an accumulation point in a Zilber field, although it has an accumulation point in $\mathbb C$.
\item
Let $\mathbb U$ be the multiplicative group of roots of unity. This has an invariant meaning in any algebraically closed field. In $\mathbb C$ any infinite subset of $\mathbb U$  has an accumulation point in $\mathbb C$  since it is a subset of the unit disc which is compact. In our more abstract situation the subset accumulates if it is infinite.

\item
Let $X$ be an infinite subset of a cyclic group $\langle \alpha \rangle$ (under multiplication).  In $\mathbb C$ we have the following three cases.

Case 1.  If $||\alpha||=1$ then $X$ has  an accumulation point on the unit circle, by compactness. If an entire function $f$ vanishes on $X$ then $f\equiv 0$.

Case 2. If $||\alpha||<1$ then $\{ \alpha^n : n\geq 0\}$ has $0$ as an accumulation point, and if $X\cap  \{ \alpha^n : n\geq 0\}$ is infinite then $X$ has $0$ as an accumulation point. Again we conclude that if $f$ is an entire function and vanishes on $X$ then $f\equiv 0$. If, however, $X\cap  \{ \alpha^n : n\geq 0\}$ is finite there is nothing we can say for a general entire function $f$, since $X$ need not have an accumulation point in $\mathbb C$ (so that Weierstrass' work \cite{alhfors} allows $X$ to be the zero set of a non identically zero entire function). Nevertheless, if we restrict $f$ to be of type (1), $\alpha$ is subject to severe constraints, relating to work of Gy\H ory and Schinzel on trinomials, see \cite{GyorySchinzel}. This is proved in Section \ref{cyclic case}, assuming Schanuel's Conjecture. The proof works uniformly for $\mathbb C$, and $\mathbb B$ (assuming (SC) for $\mathbb C$), and indeed for a much wider class of $E$-fields (the $LECS$-fields).

Case 3.  $||\alpha||>1$. This is dual to Case 2, replacing $\alpha$ by  $\frac{1}{\alpha} $ and is treated accordingly in Section \ref{cyclic case}.

\item
Let $X$ be a finite dimensional $\mathbb Q$-vector space. In $\mathbb C$ we fix a  basis and identify $X$ with $\mathbb Q^n$ where $n$ is the dimension of the space. Then we extend the definition given in 1 of Cauchy sequence to the case $n>1$ using the supremum metric. This gives us a notion of a subset of $X$ accumulating. Notice that the notion does not depend on the choice of the basis by elementary matrix theory.

\end{enumerate}

\section{Results about general infinite set of zeros}
\label{Case infinite set} 
We are going to employ some of the arguments used in \cite{DMT} for the proof of Shapiro's Conjecture from Schanuel's Conjecture.  Here we are principally interested on conditions on an infinite  set $X$ in order for it to be contained in the zero set of some exponential polynomial.  We will use these results in Section \ref{cyclic case} to deal with the case of  an infinite set of roots $X$ contained in an infinite cyclic group.
\par
 The assumptions we have through this section are 
 \begin{enumerate}
\item  
$K$ is a $LECS$-field;
\item
$X$ is an infinite subset of $K$;
\item
$tr.d._{\mathbb Q}(X)=M<\infty$;
\item $X$ is contained in the zero set of some exponential polynomial $f$ in $\mathcal E$.
 \end{enumerate}
 
Our objective is to understand the strength of these assumptions. In order to do this we proceed through various reductions, roughly following \cite{DMT}. 

\smallskip
\noindent
{\bf Reduction 1.} We apply  Ritt Factorization (Theorem \ref{ritt}) to  $f$ in $\mathcal E$, and  we use the Pigeon-hole Principle to go to the case of  $X$ infinite and  $f$  either irreducible or simple of length $2$.  We will consider the case of $f$ simple at the end of the section, and for now we assume $f$ irreducible until further notice. 

\smallskip
\noindent {\bf Reduction 2.} This involves the very first step in our proof of Shapiro's Conjecture from Schanuel's Conjecture (see Section  5 of \cite{DMT}).   Here we have the hypothesis that the transcendence degree of $X$ is finite in contrast to our previous paper on Shapiro's Conjecture where we have to prove that the set of common zeros of two exponential polynomials have finite transcendence degree.  Both in \cite{DMT} and in this paper Schanuel's Conjecture is crucial. 
\par
We also use  the same  notation and  for convenience we recall that 
\begin{enumerate}
\item[(i)]
$D=l.d._{\mathbb Q}(\mu_1,\ldots ,\mu_N)$
\item[(ii)]
$ \delta_1= tr.d._{\mathbb Q}(\lambda_1,\ldots ,\lambda_N)$
\item[(iii)]
$ \delta_2= tr.d._{\mathbb Q}(\mu_1,\ldots ,\mu_N)$
 \item[(iv)]
 $L$ is the algebraic closure of $\mathbb Q(\overline{\lambda})$
 \item[(v)] 
 $G_m^D$ is the multiplicative group variety.
 \end{enumerate}

 Let $\alpha_1, \ldots , \alpha_k\in X$ be solutions of $f(z)=0$.  An upper bound for  the linear dimension of the set $\{ \alpha_j\mu_i : 1\leq j\leq k, 1\leq i\leq N\}$ over $\mathbb Q$  is $Dk$. This is the actual dimension $D$ when  $k=1$.  We exploit the fact that Schanuel's Conjecture puts restrictions on $k$ for the above linear dimension to be $Dk$. Let $F(Y_1,\dots , Y_D)$ be the Laurent   polynomial over $L$ associated to $f$. The condition $F(Y_1,\dots , Y_D)=0$ defines an irreducible  subvariety $V$ of $G_m^D$ of dimension $D-1$ over $L$ (here we use that $f$ is irreducible in $\mathcal E$, and hence $F(Y_1,\dots , Y_D)$ is power irreducible). We get the following easy estimates 
$$tr.d._{\mathbb Q}( \overline{\mu}\alpha_1,\ldots , \overline{\mu}\alpha_k, e^{\overline{\mu}\alpha_1}  ,\ldots , e^{\overline{\mu}\alpha_k}) \leq M+\delta_2 +\delta_1+ k(D-1).$$ 
By Schanuel's Conjecture the above transcendence degree is greater or equal  than $l.d._{\mathbb Q}(\overline{\mu}\alpha_1,\ldots , \overline{\mu}\alpha_k)$, so if this dimension is $kD$ we get $$  kD\leq M+\delta_2 +\delta_1+ k(D-1).$$  Hence,  $k\leq M+\delta_2 +\delta_1$. 
\par

\smallskip
\noindent {\bf Reduction 3.}  We now appeal, as in \cite{DMT}, to work of Bombieri, Masser and Zannier \cite{bombieri} on anomalous subvarieties. 
\par
If $\alpha \in X$ we write $e^{\overline{b}\alpha}$ for the tuple $(e^{b_1\alpha}, \ldots ,e^{b_D\alpha} )$.
The variety  $V$ contains all points of the form $(e^{b_1\alpha}, \ldots ,e^{b_D\alpha} )$ for $\alpha \in X$. We consider points $\alpha_1,\ldots , \alpha_k\in X$ for $k>\delta_1+\delta_2+M$ so that among $\overline{b}\alpha_1, \ldots , \overline{b}\alpha_k$  there are nontrivial linear relations over $\mathbb Q$,  hence the linear dimension over $\mathbb Q$ of $\overline{b}\alpha_1, \ldots , \overline{b}\alpha_k$ is $<Dk$.  The point $(e^{\overline{b}\alpha_1}, \ldots ,e^{\overline{b}\alpha_k} ) $ in the $Dk$-space lies on $V^k$. 
\par
The nontrivial $\mathbb Q$-linear relations on the $b_j\alpha_r$'s induce  algebraic relations between the $e^{b_j\alpha_r}$'s. These latter relations define an algebraic subgroup $\Gamma_{\overline{\alpha}}$ of $(G_m^D)^k$ of dimension $d(\overline{\alpha})$ over $\mathbb Q$ and codimension $Dk-d(\overline{\alpha})$. The dimension $d(\overline{\alpha})$ is strictly connected to the linear dimension of the $b_j\alpha_r$'s.

The  point $(e^{\overline{b}\alpha_1}, \ldots , e^{\overline{b}\alpha_k})$ lies on $\Gamma_{\overline{\alpha}}\cap V$. Let  $W_{\overline{\alpha}}$ be  the variety of the point  $(e^{\overline{b}\alpha_1}, \ldots , e^{\overline{b}\alpha_k})$ over $L$. Now we examine the issue when $W_{\overline{\alpha}}$ is anomalous in $V^k$. See \cite{bombieri} for definitions and properties of anomalous subvarieties. 
\par
Suppose  $W_{\overline{\alpha}}$ is neither anomalous nor of dimension $0$.  Then, in particular,  $\dim W_{\overline{\alpha}}\leq \dim (V^k)-  \mbox{ codim } (\Gamma_{\overline{\alpha}}), $ i.e. $\dim W_{\overline{\alpha}}\leq  k(D-1) - (Dk-d(\overline{\alpha}))= d(\overline{\alpha})-k$. Schanuel's Conjecture implies 
 $d(\overline{\alpha})\leq d(\overline{\alpha})-k+\delta_1+M+\delta_2$, hence $k\leq \delta_1+ M + \delta_2$. So, if $k> \delta_1+ M + \delta_2$ then $W_{\overline{\alpha}}$  is either anomalous or of dimension $0$.

\smallskip
\noindent {\bf Reduction 4.} Suppose such a $W_{\overline{\alpha}}$ has dimension $0$. Then all the $e^{b_j\alpha_r}$'s are algebraic over $L$. So, $tr.d._{\mathbb Q}(e^{\overline{b}\alpha_1}, \ldots , e^{\overline{b}\alpha_k})\leq \delta_1$. Moreover, $tr.d._{\mathbb Q}(\overline{b}\alpha_1, \ldots , \overline{b}\alpha_k)\leq \delta_2+M.$ So, by Schanuel's Conjecture, 
\begin{equation}
\label{bounddalpha}
 d(\overline{\alpha})\leq \delta_1+\delta_2+M. 
\end{equation}

\smallskip
\noindent {\bf Reduction 5.} We worked locally with $\overline{\alpha}$ and now we have to work independently of $\overline{\alpha}$. Notice that under permutations of the $\overline{\alpha}$ the two properties, that $W_{\overline{\alpha}}$ has dimension $0$, or  $W_{\overline{\alpha}}$ is anomalous, are invariant. Suppose $l.d._{\mathbb Q}(X)$ is infinite. Choose an infinite independent subset $X_1$ of $X$. Let $k>\delta_1+\delta_2+M$, and $\overline{\alpha}=(\alpha_1,\ldots ,\alpha_k)$, with $\alpha_1,\ldots ,\alpha_k$ distinct elements of $X_1$, hence linearly independent over $\mathbb Q$. Then  $W_{\overline{\alpha}}$ cannot have dimension $0$ since $d(\overline{\alpha})\geq k$ and (\ref{bounddalpha}) holds. Thus, for any $k$-element subset $\{ \alpha_1,\ldots ,\alpha_k\}$ of $X_1$, $W_{\overline{\beta}}$ is anomalous in $V^k$ for any permutation $\overline{\beta}$ of $( \alpha_1,\ldots ,\alpha_k)$. We thin $X$ to $X_1$ and we work with $X_1$.

\smallskip
\noindent {\bf Reduction 6.} We work with $X_1$ which we still call $X$. Let $k$ be  the minimal such that for any $k+1$ elements $\eta_1, \ldots , \eta_{k+1}$ of $X$, the variety of the point $(e^{\overline{b}\eta_1}, \ldots , e^{\overline{b}\eta_{k+1}})$ is anomalous in $V^{k+1}$. From \cite{bombieri} it follows that there is a finite collection $\Phi$ of proper tori $H_1,\ldots ,H_t$ is $G_m^{D(k+1)}$ such that each maximal anomalous subvariety of $V^{k+1}$ is a component of the intersection of $V^{k+1}$ with a coset of one of the $H_i$. We now discard much information. For each $H_i$ we pick one of the multiplicative conditions defining it. These define a finite set $\{J_1,\ldots, J_t\}$ of codimension $1$ subgroups so that each anomalous subvariety is contained in one of them. Now, exactly as in \cite{DMT}, by using Ramsey's Theorem and Schanuel's Conjecture yet again, we get an infinite $X_2\subseteq X$ so that the $\mathbb Q$-linear dimension of $X_2$ is finite.

\smallskip
\noindent {\bf Reduction 7.} Without loss of generality we can assume that  the set $X$ of solutions of the exponential polynomial $f$ is infinite and of finite linear dimension over $\mathbb Q$. Each $\alpha$ in $X$ gives rise to a solution $(e^{\mu_1\alpha}, \dots , e^{\mu_N\alpha})$ of the linear equation $\lambda_1Y_1+\ldots +\lambda_NY_N=0$. We change the latter to the equation 

\begin{equation}
\label{linear equation}
\biggl(\frac{-\lambda_1}{\lambda_N}\biggr)Z_1+ \ldots + \biggl(\frac{-\lambda_{N-1}}{\lambda_N}\biggr) Z_{N-1}=1
\end{equation}
where $Z_j$ stands for $\frac{Y_{j}}{Y_N}$. Note that $(e^{(\mu_1-\mu_N)\alpha}, \dots , e^{(\mu_{N-1}-\mu_N)\alpha})$ is a solution of equation (\ref{linear equation}).   In \cite[Lemma 5.6]{DMT} we observed that distinct $\alpha$'s give distinct roots of (\ref{linear equation}), unless $f$ is a simple polynomial. The finite dimensionality of $X$ implies that the multiplicative group generated by  $e^{(\mu_{j}-\mu_N)\alpha}$'s for $\alpha\in X$ and $j=1,\ldots N-1$ has finite rank. We can then apply a basic result on solving linear equations over a finite rank multiplicative group,  due to Evertse, Schlickewei and Schmidt in  \cite{ESS}. From this result it follows  that only finitely many solution of (\ref{linear equation})  of the form $(e^{(\mu_1-\mu_N)\alpha}, \dots , e^{(\mu_{N-1}-\mu_N)\alpha})$, for $\alpha \in X$,  are non-degenerate. We then thin $X$ again to an infinite set, which we still call $X$, generating infinitely many degenerate solutions of (\ref{linear equation}).  For any proper subset $I$ of $\{ 1,\ldots , N-1\}$ with $|I|>1$ let $$f_I =\sum_{j\in I}\biggl(\frac{-\lambda_j}{\lambda_N}\biggr) e^{(\mu_j-\mu_N)z}.$$
 By repeated applications of the Pigeon-hole Principle and  the Evertse, Schlickewei and Schmidt result in  \cite{ESS} we construct a finite chain of subsets $I_j$ of $\{  1,\ldots ,N-1\} $  so that $f_{I_j}$ has infinitely many common solutions with $f$. For cardinality reasons we have to reach an $I_{j_0}$ of cardinality $2$ whose corresponding polynomial is simple of length $2$.  By our result \cite[Theorem 5.7]{DMT} we get that $f$ divides a simple polynomial, and so $f$ is necessarily simple.  
\par
Recall that in Reduction 1 we postponed the discussion for simple polynomials, and we did the reductions for irreducible polynomials in $\mathcal E$. Now we have reached the conclusion that only a simple polynomial can satisfy the assumptions 1 - 4 of this section.  We have then proved the following 

\begin{thm}
\label{main}
(SC) Let $f$ be an exponential polynomial in $\mathcal E$.   If  $Z(f)$ contains an  infinite set $X$ of finite transcendence degree then $f$ is divisible by a simple polynomial. Every infinite subset of the zero set of an irreducible polynomial has infinite transcendence degree (hence infinite linear dimension). 
\end{thm}

\begin{osserva}
{\rm 
 Notice that the above theorem and Lemma \ref{simplepoly} imply that if an exponential polynomial vanishes on an infinite set $X$ either of finite transcendence degree or of finite dimension  then, assuming Schanuel's Conjecture in the first case and unconditionally in the second, the linear dimension of $X$ is at most $2$.  
The above results also imply  that  the zero set of an irreducible exponential polynomial cannot contain an infinite set  of finite linear dimension over $\mathbb Q$. This does not depend on Schanuel's Conjecture. So there is no irreducible exponential polynomial $f$ such that $Z(f)$ contains an infinite set of algebraic numbers. If $Z(g)$ contains an infinite set of finite dimension  and $g$ is simple of length $2$ then  this dimension has to be at most $2$.}
\end{osserva}

\subsection{The case of subsets of $\mathbb Q$}
\label{rationalsolutions}

An immediate consequence of the above results is that if $f$ is an exponential polynomial whose zero set contains an infinite subset $X$ of rationals then $f$ is divisible by a simple polynomial. We can then assume without loss of generality that $f$ is simple. We now show that  the elements in $X$ have bounded denominators. Indeed,  factor $f$ into simple polynomials $f_1, \ldots, f_k$ of length $2$, and let $X_i=X\cap Z(f_i)$. For any $i$,   $X_i\subseteq Z(f_i)\cap \mathbb Q$ and  moreover,    $X_i\subseteq A_{f_i}+ \mathbb ZB_{f_i}$ in the notation of Lemma \ref{simplepoly}. If $q_1,q_2 \in X_i$ then $q_1-q_2 = (k-h)B_{f_i}$ for some $k,h\in \mathbb Z$. Hence, $B_{f_i}\in \mathbb Q$, and so also $A_{f_i}\in \mathbb Q$.  This implies that there is $N_i\in \mathbb N$ such that $X_i\subseteq \frac{1}{N_i}\mathbb Z$. So, $X$ has bounded denominators since $X=X_1\cup \ldots \cup X_k$.
\par
There is a  recently published paper by Gunaydin \cite{gun} about solving in the rationals exponential  polynomials in many variables $\overline{X} =(X_1, \ldots , X_t)$ over  $\mathbb C$ of the form 
\begin{equation}
\label{polyGunaydin}
\sum_{i=1}^{s}P_i(\overline{X})e^{((\overline{X}\cdot \overline{\alpha_i}))}
\end{equation}
where $P_i(\overline{X})\in \mathbb C[\overline{X}]$  and $\overline{\alpha_i}\in \mathbb C^t.$ 
His result is 
\begin{thm}\cite{gun} 
\label{gun} 
Given $P_1, \ldots, P_s\in \mathbb C[\overline{X}]$ and $\overline{\alpha_1}, \ldots , \overline{\alpha_s}\in \mathbb C^t$, there is $N\in \mathbb N^{>0}$ such that if $\overline{q}\in \mathbb Q^t$ is a non degenerate solution of $$\sum_{i=1}^{s}P_i(\overline{X})e^{((\overline{X}\cdot \overline{\alpha_i}))}=0$$ then $\overline{q}\in(\frac{1}{N}\mathbb Z)^t$.
\end{thm}

He makes no reference to any other exponential field but we have verified that his results hold for exponential polynomials as in (\ref{polyGunaydin}) over a $LEC$-field.

His conclusion for $f(z)$ in a single variable implies that the rationals in the zero set of the function do not {\it accumulate}. So we have an algebraic proof of a result  proved for $\mathbb C$ via analytic methods.

\noindent {\bf Open problem.}  It is a natural question to ask if the result of Theorem \ref{gun} can be extended  to an arbitary exponential  polynomial with  iterations of exponentiation over  a $LEC$-field or more generally for a Zilber field.

\smallskip
If Zilber's Conjecture is true there is no exponential polynomial in this general sense  vanishing on a set of rationals which has an accumulation point.  
Note that there are infinite subsets of $\mathbb Q$ with unbounded denominators, without accumulation points, e.g. $$X=\biggl\{ \frac{m_n}{p_n} : p_n \mbox{ is $n$th prime, $m_n\in \mathbb N$,  and $n<\frac{m_n}{p_n}$} < n+1\biggr\}. $$
It seems inconceivable to us that there is a non trivial exponential polynomial vanishing on $X$. The classic  Weierstrass Theorem (see \cite{alhfors}) provides a nontrivial entire function vanishing on $X$, but the usual proof  gives a complicated infinitary definition. 

In the special case of exponential polynomials over $\mathbb C$ of the form $$  \sum_{i=1}^{s} \lambda_ie^{\mu_ie^{2\pi iz}}$$ Theorem \ref{gun} is true because otherwise the polynomial $$  \sum_{i=1}^{s}  \lambda_ie^{\mu_iw}$$ would have infinitely many roots of unity as solutions. We prove in the next section that this cannot happen unless the polynomial is identically zero.

\section{Case of roots of unity} 

Let $\mathbb U$ denote the set of roots of unity. We know that an entire function over $\mathbb C$ cannot have infinitely many roots of unity as zeros unless it is identically zero. Let $K$ be a $LECS$-field.

\begin{thm}
\label{rootsofunity}
(SC) If $f(z)\in \mathcal E$ over $K$ vanishes on an infinite subset $X$ of $\Bbb U,$ then $f(z)$ is the zero polynomial in $\mathcal E$. 
\end{thm}

\noindent \textbf{Proof:} Without loss of generality $f$ is either irreducible or simple of length $2$.
\par
\noindent
{\it Case 1:}  Let $f(z) = \lambda_1 e^{\mu_1 z} + \ldots + \lambda_N e^{\mu_N z}$ be irreducible. Clearly, $tr.d._{\mathbb Q}(X)=0$. The same arguments,  including Schanuel's Conjecture, used in the proof of Theorem \ref{main} give  the finite dimensionality of $X$ over $\mathbb Q$.  Let $\{ \theta_1, \ldots, \theta_k \}\subseteq  X$ be  a   $\Bbb Q$-basis of $X$.  Any other solution of $f$ in $X$ belongs to the finitely generated field $F=\mathbb Q(\theta_1, \ldots, \theta_k )$.  It is a very well known fact that any finitely generated field $F$ of characteristic $0$  cannot contain infinitely many roots of unity, see \cite{Bourbaki}.  So we get a contradiction.  

\noindent
{\it Case 2:}  If $f$ is a simple polynomial of length $2$ then Lemma  \ref{simplepoly} gives unconditionally that the  linear dimension of $X$ is finite, forcing $X$ to be finite,  and we get again a contradiction. \hfill \qed

\subsection{Torsion points of elliptic curves}

We generalize the above argument to the case of coordinates of torsion points of an elliptic curve $E$ over $\mathbb Q$.  We need the following well known result, see \cite{Bourbaki}.

\begin{thm}
Let $F$ be a characteristic $0$ field. If $F$ is  finitely generated  then any subfield of $F$ is also finitely generated.

\end{thm}

In Merel's paper \cite{Merel} the following major result is proved. 

\begin{thm}
\label{Merel}
For all $d\in \mathbb Z$, $d\geq 1$ there exists a constant $n(d)\geq 0$ such that for all elliptic curves  $E$ over a number field $K$ with $[K:\mathbb Q]=d$ then any torsion point of $E(K)$ is of order less than $n(d)$. 

\end{thm}

Recall that for each $n$ the group of $n$-torsion points of an elliptic curve $E$ is isomorphic to $\mathbb Z/n\mathbb Z\oplus \mathbb Z/n\mathbb Z$. 

\smallskip
\noindent
{\bf Setting.} We consider the affine part of the curve given as usual by an equation quadratic in $x$ and cubic in $y$. Define $\pi_1$ as the projection to the $x$-coordinate, and $\pi_2$ as the  projection to the $y$-coordinate. 

Let $E(\mathbb Q)_{tors}$ be the torsion points of $E$ in $\mathbb Q^{alg}$, and $X= \pi_1 (E(\mathbb Q)_{tors})\cup  \pi_2 (E(\mathbb Q)_{tors})$. 

\begin{thm}
\label{elliptic}(SC) 
Let $f$ be an exponential polynomial in $\mathcal E$. Then $X\cap Z(f)$  is finite. 

\end{thm}

\noindent \textbf{Proof:} Suppose $X\cap Z(f)$ is infinite. By earlier remarks we can assume $f$ is simple, and as in the proof of Theorem \ref{rootsofunity} we can assume that we have an infinite subset $X_1\subseteq X\cap Z(f)$ of finite dimension over $\mathbb Q$. Hence, there is a number field $F$ of dimension $d$ such that $X_1\subseteq F$. So by Merel's result any torsion point of $E$ over $F$ has order $\leq n(d)$.  So, there are only finitely many of them and then we get a contradiction.  \hfill \qed

\medskip
A natural generalization of Theorem \ref{elliptic} is stated in the following 

\smallskip

\noindent 
{\bf Open Problem:}  Is the intersection of an infinite set of solutions of a non zero exponential polynomial with the coordinates of the torsion points of an abelian variety over $\mathbb Q$  finite?

\smallskip
The crucial obstruction is that as far as we know there is not an analogous result to Merel's for abelian varieties.

\section{Case of an infinite cyclic group}
\label{cyclic case}
In this section we examine the case of  an infinite $X\subseteq Z(f)$ such that $X\subseteq \langle \alpha \rangle$, where $\alpha \in K$. Clearly, the transcendence degree of $X$ is finite,  and by Theorem \ref{main} without loss of generality we can assume $f$ is  simple (modulo Schanuel's Conjecture). We are going to show that some power of $\alpha$ is in $\mathbb Z$. Without loss of generality by factorization $f$ can be chosen of length $2$. 
\par
In this section we will often use the following basic thinning and reduction argument for the set $X$ of solutions using the Euclidean division. 

\smallskip

\noindent 
{\bf Euclidean reduction.}   Suppose $\alpha^r \in X$ for some $r\in \mathbb N$ (in case $r$ is negative we work with $\frac{1}{\alpha}$).  Let $s_0<r$ such that there are infinitely many $m$'s with $\alpha^{rm+s_0}\in X$. Via a change of variable we work with the  polynomial   $g(z)=f(\alpha^{s_0}z)$ which vanishes on  an infinite subset of $\langle \alpha^r \rangle$. If $f$ is simple of length $2$ then also $g$ is simple of length $2$. Hence by Lemma \ref{simplepoly}, $Z(g)$ is the translate $A_g+ \mathbb ZB_g$, where $A_g$ and $B_g$ are in $K$,  and  in general different from $A_f$ and $B_f$, respectively. In what follows it will not make any difference if we work with either $f$ or $g$. The infinite set $X$ of solutions of $f$ contains a translate $X'$ of an infinite subset of $Z(g)$. We will not make any distinction between the two, and we will still use the notation $X$ for $X'$.

\begin{lemma}
\label{integersolutions}
Let $\alpha \in K^{\times}$ and $\alpha \not\in \mathbb U$. Suppose that the set $X\cap \langle \alpha \rangle$ is infinite. 
If for some $r\not=0$ in $\mathbb Z$, $\alpha^r\in \mathbb Q$ then $\alpha^r\in \mathbb Z$ or $\frac{1}{\alpha^r}\in \mathbb Z$. 

\end{lemma}

\noindent \textbf{Proof:} Suppose $\alpha^r$ is rational and belongs to $  X\cap \langle \alpha \rangle$.   Via the Euclidean reduction we can assume that there are infinitely many $m$'s with $\alpha^{rm}$ rational solutions of $g$.  Choose two of them, say $\alpha^{rm_1}$ and $\alpha^{rm_2}$, so  $\alpha^{rm_1}=A_g+kB_g$ and  $\alpha^{rm_2}=A_g+hB_g$, for some $k,h\in \mathbb Z$. Clearly, $A_g$ and $B_g$ are rationals.  If there is no upper bound on the $m$'s then we get a contradiction since the $\alpha^{rm}$'s are contained in $\frac{1}{N}\mathbb Z$ for some $N$ (see Section \ref{rationalsolutions}). If there is no lower  bound on the $m$'s then it is sufficient to work with $\frac{1}{\alpha}$ as a generator. It follows that either $\alpha^r\in \mathbb Z$ or $\frac{1}{\alpha^r}\in \mathbb Z$.  \qed

\begin{thm}
\label{cyclicgroupcase}
Let $f\in \mathcal E$. Suppose that $X\subseteq Z(f)$, $X$ is contained in the infinite cyclic group generated by $\alpha$, and $X$ is infinite. Then  $f$ is identically zero unless $\alpha^r\in \mathbb Z$ for some $r\in \mathbb Z$. 
\end{thm}

\noindent \textbf{Proof:}  As already observed $f$ is simple of length $2$, and so $Z(f)=A_f + \mathbb ZB_f$. Hence $X$ has linear dimension over $\mathbb Q$ less or equal than $2$.  In particular, $\alpha $ is  a root of infinitely many trinomials over $\mathbb Q$, and so it is algebraic with minimum polynomial $p(x)$ dividing infinitely many trinomials.   By work of Gy\"ory and Schinzel in \cite{GyorySchinzel}, there exists a polynomial $q(x)\in \mathbb Q[x]$ of degree $\leq 2$ such that $p(x)$ divides  $q(x^r)$ for some $r$.  
\par
\noindent If $q(x)$ is linear then $\alpha^r$ is rational and we have finished thanks to  Lemma \ref{integersolutions}.  
\par
\noindent If $q(x)$ is quadratic  then  $\alpha^r \in \mathbb Q(\sqrt{d})$ for some $d\in \mathbb Q$.  By  Euclidean reduction we  reduce to the case of a simple polynomial $g$ of length $2$ for which the corresponding  $A_g$ and $B_g$ are in  $\mathbb Q(\sqrt{d})$. 
Using again  Euclidean reduction we can assume without loss of generality that  $\alpha \in \mathbb Q(\sqrt{d})$. The polynomial $g$ may have changed but we will continue to refer to it as $g$. Note that if $\mathbb Q(\sqrt{d})=\mathbb Q $ then $\alpha \in \mathbb Q$, and so by Lemma \ref{integersolutions} either $\alpha \in \mathbb Z$ or $\frac{1}{\alpha} \in \mathbb Z$. So, we can assume $\mathbb Q(\sqrt{d})\not=\mathbb Q $

We will also drop the subscript $g$ from $A_g$ and $B_g$ in the rest of the proof since no confusion can arise. 

\smallskip
\noindent
{\it Claim.} Either $\alpha$ or  $1/\alpha$ is  an algebraic integer in $\mathbb Q(\sqrt{d})$.

\noindent
{\it Proof of the Claim.}Suppose $\alpha$ is not an algebraic integer in $\mathbb Q(\sqrt{d})$. By a basic result in algebraic number theory (e.g. see \cite{cassels}) there is  a valuation $v$ on $\mathbb Q(\sqrt{d})$ such that $v(\alpha)<0$. If for infinitely many positive $m\in \mathbb Z$, $\alpha^m\in X$ then  there is no lower bound on the valuations of  elements of $X$. We get a contradiction since for all $m$, $\alpha^m=A+k_mB$ for some $k_m\in \mathbb Z$, and $v(A+k_mB)\geq \min \{ v(A), v(B)\}$. So the valuations of elements of $X$ have to be bounded below.  If the infinitely many integers  $m$ such that  $\alpha^m\in X$ are negative  then apply the same argument to $\frac{1}{\alpha}$. 

\smallskip
Let $\alpha^m\in X$ so that $\alpha^m=A+k_mB$ for some $k_m\in \mathbb Z$, and let $\sigma$ be a generator of the Galois group of $\mathbb Q(\sqrt{d})$ over $\mathbb Q$. Then $\sigma (\alpha)^m= \sigma (A)+k_m\cdot \sigma(B)$. The norm function is defined as  $Nm(\alpha)=\alpha \cdot \sigma (\alpha) $.   So, $$Nm(\alpha)^m= Nm(A)+Tr(A\cdot \sigma (B))k_m+Nm(B)k_m^2,$$ 
where $Tr$  denotes  the trace function. The polynomial $Nm(A)+Tr(A\cdot \sigma (B))x+Nm(B)x^2$ is over $\mathbb Q$.  By Euclidean reduction with $r=3$ and some $s_0$ with $0\leq s_0<3$, the equation 
\begin{equation}
\label{normequation}
Nm(\alpha^{s_0})y^3= Nm(A)+Tr(A\cdot \sigma (B))x+Nm(B)x^2
\end{equation}
has infinitely many integer solutions of the form $(Nm(\alpha^m), k_m )$. Notice that according to the sign of the infinitely many $m$'s we work either  with $\alpha$ or $\frac{1}{\alpha}$. If the polynomial $$P(x)=Nm(A)+Tr(A\cdot \sigma (B))x+Nm(B)x^2$$ has nonzero discriminant then (\ref{normequation}) defines (the affine part of) an elliptic curve  over $\mathbb Q$, and by Siegel's Theorem (see \cite{bombieri1}) we get a contradiction. 
\par
It remains to consider the case when the discriminant of $P(x)$ is zero, i.e. $ (Tr(A\cdot \sigma (B)))^2-4Nm(A\cdot B)=0$. In this case we have $$P(x)= Nm(B)\biggl(x+\frac{Tr(A\cdot \sigma (B))}{2Nm(B)}\biggr)^2$$ where $-\frac{Tr(A\cdot \sigma (B))}{2Nm(B)}$ is the multiple root, and equation (\ref{normequation}) becomes 
\begin{equation}
\label{reducednormequation}
Nm(\alpha^{s_0})y^3= Nm(B)\biggl(x+\frac{Tr(A\cdot \sigma (B))}{2Nm(B)}\biggr)^2.
\end{equation}
Equation (\ref{reducednormequation}) has infinitely many rational solutions of the form 
\begin{equation}
\label{solutions}
\biggl(Nm(\alpha)^m, k_m+\frac{Tr(A\cdot \sigma (B))}{2Nm(B)}\biggr) 
\end{equation}
as $m$ varies in $\mathbb Z$. By the change of  variable $x\mapsto x+\frac{Tr(A\cdot \sigma (B))}{2Nm(B)}$ and dividing by $Nm(\alpha^{s_0})$, we transform equation
(\ref{reducednormequation}) to one of the following form $y^3=cx^2$, where $c=\frac{Nm(B)}{Nm(\alpha^{s_0})} \in \mathbb Q^{\times}$. The equation defines a rational curve which we now parametrize. 

\smallskip

\noindent
{\it Claim.} The rational solutions of $y^3=cx^2$ are of the following form
\begin{equation}\label{parametrizzazione}
\left\{ \begin{array}{ll} 
x =  &  \theta^3c^{-5}  \\
y =  & \theta^2c^{-3} 
\end{array} \right.\\
\end{equation}
with $\theta \in \mathbb Q$.

\noindent
{\it Proof of the Claim:} 
Consider the $p$-adic valuation $v_p$ for some $p$.  Suppose  $y^3=cx^2$ with $x,y\in \mathbb Q$ then  we have 
\begin{enumerate}
\item
$3v_p(y)=2v_p(x)+v_p(c), \mbox{ and so } v_p(y)=v_p\Bigl( \bigl(\frac{x}{y}\bigr)^2c\Bigr)$
\item
$15v_p(y)=10v_p(x)+v_p(c), \mbox{ and so } v_p(x)= v_p\Bigl( \bigl(\frac{y^5}{x^2}\bigr)^3c^{-5}\Bigr)$. 
\end{enumerate}

Hence,  $x=\theta^3c^{-5}$  and $y= \xi^2c$, for some $\theta,\xi$, and so $\xi^6c^3= \theta^6c^{-9}$. Therefore, $(\xi /\theta)^6 =c^{-12}$, and this implies $\xi /\theta=\pm c^{-2}$, i.e. $\xi = \pm \theta c^{-2}$. Notice that in the equation $y^3=cx^2$ the variable $x$ occurs in  even power while the variable $y$ occurs in  odd power, so  we get the following parametrization of the curve $y^3=cx^2$
\begin{equation}\label{parametrizzazione}
\left\{ \begin{array}{ll} 
x =  & \pm \  \theta^3c^{-5}  \\
y =    & \theta^2c^{-3} 
\end{array} \right.\\
\end{equation}

We have then obtained that all the rational solutions of $y^3=cx^2$ are of  the form in (\ref{parametrizzazione}) for $\theta $ rational. 
From (\ref{solutions}) and (\ref{parametrizzazione}) it follows that for infinitely many $m\in \mathbb N$ 

\begin{equation}
\label{pair}
\biggl(Nm(\alpha)^m, k_m+\frac{Tr(A\cdot \sigma (B))}{2Nm(B)}\biggr) = ( \theta_m^2c^{-3} , \theta_m^3c^{-5}),
\end{equation}
 for $\theta_m\in \mathbb Q$. Hence $k_m = \theta_m^3c^{-5}-\frac{Tr(A\cdot \sigma (B))}{2Nm(B)}$. 
\par 
 By the Euclidean reduction  there are infinitely many $m\in \mathbb N$ such that $(\alpha^m)^3\alpha^{s_0}= A+k_mB$. Now we consider the new curve  in $y$ and $\theta$, 
\begin{equation}
\label{newequation}
y^3\alpha^{s_0}= \theta^3c^{-5}B+ \biggl(A-\frac{Tr(A\cdot \sigma (B))}{2Nm(B)}B\biggr).
\end{equation}
We distinguish now two cases. 
\smallskip

\noindent
Case 1: $A-\frac{Tr(A\cdot \sigma (B))}{2Nm(B)}B=0$. Let $m,m'\in \mathbb Z$ such that  $\alpha^m\not= \alpha^{m'}$ and both belong to $X$, then $$\alpha^{m-m'}=\biggl(\frac{\theta_m}{\theta_{m'}}\biggr)^3$$ and by Lemma \ref{integersolutions} we complete the proof in  this case.

\smallskip 
\noindent 
Case 2: $A-\frac{Tr(A\cdot \sigma (B))}{2Nm(B)}B\not=0$. Dividing on both sides of equation (\ref{newequation}) by $\alpha^{s_0}$   we obtain a new elliptic curve defined on the number field $\mathbb Q(\sqrt{d})$.  In this case we apply a generalization of Siegel's theorem (see \cite{bombieri1}) to the new elliptic curve 
defined on $\mathbb Q(\sqrt{d})$. Fix a prime ideal $\mathcal P$ in $\mathbb Q(\sqrt{d})$ then by (\ref{pair}) we have that  $$v_{\mathcal P}(\theta_m) \geq \frac{1}{3}\biggl(5v_{\mathcal P}(c)+\min \biggl\{0,v_{\mathcal P}\biggl(\frac{Tr(A\cdot \sigma (B))}{2Nm(B)}\biggr)\biggr\}\biggr).$$  For all but finitely many prime ideals  $\mathcal P$ we have both $v_{\mathcal P}(c)$ and $v_{\mathcal P}\biggl(\frac{Tr(A\cdot \sigma (B))}{2Nm(B)}\biggr)\geq 0$. So there is a finite set  $\mathcal S$ of prime ideals $\mathcal P$ in $\mathbb Q(\sqrt{d})$  so that all $v_{\mathcal P}(\theta_m)$ are $\mathcal S$-integer. Recall also that $\alpha^{s_0} $ is an integer in $\mathbb Q(\sqrt{d})$, so we get a contradiction with Siegel's theorem.
\qed

\smallskip

We are now able to characterize those infinite subsets of an infinite cyclic group which may occur as zero-set of an exponential polynomial in $\mathcal E$. In particular, this polynomial is divisible by  a  simple polynomial. 

\begin{corollario}
The following characterization holds $$ \{ \alpha \in K: 
|\langle \alpha \rangle \cap Z(f)|=\infty  \mbox{ {\rm for some} } f\in \mathcal E\} =\{ \alpha \in K : \alpha^r \in \mathbb Z \mbox{  {\rm for some} } r\in \mathbb Z - \{ 0\} \}.$$
\end{corollario}

\noindent
{\bf Proof:} $(\subseteq )$  This inclusion follows from previous theorem. 

\noindent
$(\supseteq  )$  Clearly $\alpha$ is algebraic over $\mathbb Q$. Suppose $r\in \mathbb Z$ and $r>0 $ then for any $s\geq r$ $\alpha^s\in \mathbb Z$. Let $j_0<r$ and consider the infinitely many integers of the form $s=rq+j_0$ as $q$ varies in $\mathbb N$. Then the polynomial $f(z)= 1-e^{\mu z}$, where $\mu= \alpha^{j_0}2i\pi$,  has infinitely many roots in $\langle \alpha \rangle$.  If $r<0$ then use the same argument for $\frac{1}{\alpha}$. \qed

\bigskip
\noindent
{\bf Acknowledgements.} This research was conducted at ICMS in Edinburgh. The authors thank ICMS  for funding, and the ICMS staff for their generous practical support. This research is part of  the project FIRB 2010, {\it Nuovi sviluppi nella Teoria dei Modelli dell'esponenziazione}.

\end{document}